\documentclass[12pt,a4paper]{article}
\usepackage{amsmath,amssymb,amsthm,array,graphicx}
\usepackage{bbm} 
\begin{document}
\setlength{\parindent}{0pt}
\setlength{\parskip}{6pt plus 1pt}
\newcommand*{\E}{\mathbb{E}}
\renewcommand*{\P}{\mathbb{P}}
\newcommand*{\cov}{\mathrm{cov}}
\newcommand*{\Var}{\mathbb{V}\mathrm{ar}}
\newcommand*{\ind}{\mathbbm{1}}
\renewcommand{\arraystretch}{1.5}
\providecommand{\keywords}[1]{\textit{Keywords.} #1}
\providecommand{\subjclass}[1]{\textit{2010 Mathematics Subject Classification.} #1}
\newtheorem{theorem}{Theorem}
\newtheorem*{conj}{Conjecture}
\theoremstyle{remark}
\newtheorem*{remark}{\noindent Remark}

\date{27 June 2024}
\title{Three classical probability problems: the Hungarian roulette}
\author{Tam\'as F. M\'ori, G\'abor J. Sz\'ekely\\
\small{HUN-REN Alfr\'ed R\'enyi Institute of Mathematics,}\\ 
\small{Budapest, Hungary}}

\maketitle
\begin{abstract}
In this paper three unrelated problems will be discussed. What connects them is 
the rich methodology of classical probability theory. In the first two problems we 
have a complete answer to the problems raised; in the third case, what we call 
the Hungarian roulette problem, we only have a conjecture with heuristic justification. 
\end{abstract}

\keywords{
Normal distribution, exponential distribution, independence, distance 
correlation, asymptotic logarithmic periodicity of probabilities, 
merging of probabilities
}

\subjclass{60E10, 60F05, 62E10}

section{Introduction}
This paper is devoted to three unrelated problems from the field of classical 
probability theory. They are interesting individually.

The first problem leads to a joint characterization of two families of probability  
distributions: the exponential and the normal one. No doubt, they belong to the 
most important families of distributions. Both have several characterizations in 
the literature, see e.g \cite{KLR,GK,Bryc}, but a joint characterization has not 
been presented until now, as far as we know.
 
In the second problem we ask if the independence of $X-X'$ and $Y-Y'$ where 
$(X, Y)$ and $(X', Y')$ are iid, implies the independence of $X$ and $Y$.
Ex.~15.2 on p.208  of \cite{SzR} shows that $\cov\left(|X-X'|,\,|Y-Y'|\right)=0$ 
does not imply the independence of $X$ and $Y$. The joint density
\[
p(x,y)=\left(1/4-q(x)q(y)\right)\ind_{[-1,1]^2}(x,y)
\]
where
\[
q(x)=-(c/2)\ind_{[-1,0]}(x)+(1/2)\ind_{(0,c]}(x)
\]
with $c=\sqrt{2}-1$ is a counterexample. (Here and in the sequel $\ind_S$ and 
$\ind(A)$ stand for the indicator of the set $S$ and that of the random event $A$, 
resp.)

The third problem is about a game we call the Hungarian roulette. This simple
model is closely connected to classical allocation problems \cite{KSC}. Though they 
have a rich history, the question we will address seems to be missing from the 
literature.

\section{Joint characterization of the exponential and the normal distribution}

Let the random variables $X$ and $Y$ be independent and identically distributed 
with probability density function $f$ that is positive over a bounded or infinite open 
interval $I$ (and $0$ outside). 
Determine all possible distributions of $X$ for which the conditional density of 
$U := X+Y$ given $X-Y=v$ is of the form $Af(Au +B|v|)$ for some real numbers 
$A>0$ and $B$.

\begin{theorem}\label{thm1}
The probability density function $f$ satisfies the conditions above if and only if the
distribution of $X$ is either exponential, or negative exponential (that is, the distribution 
of $-X$ is exponential), or normal.
\end{theorem}

This is a novel joint characterization of the exponential
and the normal (Gaussian) distributions.

\proof
It is easy to check that the ``if'' part is true (see \eqref{feqn1} below).
In order to prove the ``only if'' direction let us compute the conditional density of
$U=X+Y$ given $X-Y=v$. 

The joint density of $U=X+Y$ and $V=X-Y$ is
\[
\frac{1}{2}\,f\Big(\frac{u+v}{2}\Big)f\Big(\frac{u-v}{2}\Big),
\]
and the density of $V$ is
\[
g(v)=\frac{1}{2}\int_{-\infty}^{\infty}f\Big(\frac{u+v}{2}\Big)f
\Big(\frac{u-v}{2}\Big)\,\mathrm du.
\]
Clearly, $g(v)>0$ if and only if $|v|<|I|$, the length of the interval $I$.

From this we obtain the following functional equation.
\begin{equation}\label{feqn1}
f\Big(\frac{u+v}{2}\Big)f\Big(\frac{u-v}{2}\Big)=2Af(Au + B|v|)g(v).
\end{equation}
This holds by supposition if $|v|<|I|$. In the opposite case at least one of
the arguments $\frac{u+v}{2}$ and $\frac{u-v}{2}$ falls outside of $I$, thus
both sides of \eqref{feqn1} are zero. Thus \eqref{feqn1} holds for all real 
numbers $u$ and $v$. We can suppose that $v\ge 0$.

Put $s=\dfrac{u+v}{2}$, $t=\dfrac{u-v}{2}$. Then \eqref{feqn1} is equivalent to
\begin{equation}\label{feqn2}
f(s)f(t)=2Af\!\big((A+B)s+(A-B)t\big)g(s-t),
\end{equation}
for all real numbers $s,t$. There are two cases.

\textbf{Case (i)}:  $A=|B|$. Suppose $A=B$. (A similar argument works if $A=-B$.)
Then \eqref{feqn2} reduces to
\begin{equation}\label{feqn3}
f(s)f(t)=2Af(2As)g(s-t).
\end{equation}
If $A=1/2$, then \eqref{feqn3} implies that both $f$ and $g$ are constant. Thus 
random variables $X$ and $Y$ are uniformly distributed on the bounded interval 
$I$. But in this case $g$ is not constant, as it has a triangular shape. This
means that $A\ne 1/2$.

If $A>1/2$ and $I=(a,b)$ with a positive and finite $b$ we have $2As>s$ for all
$s\in(0,b)$. Thus there exists an $s$ such that $a<s<b<2As$, and $f(2As)=0$. This, 
together with \eqref{feqn3} implies $f(s)=0$, which is excluded. Therefore, if 
$b>0$ then $b=+\infty$. Similarly, if $a<0$, then $a=-\infty$, and if $a\ge 0$, 
then $a=0$. We shall see that $I=(-\infty,+\infty)$ cannot hold, therefore either
$I=(-\infty,0)$ or $I=(0,+\infty)$. The same argument applies if $A<1/2$.

From \eqref{feqn3} we get
\begin{equation}\label{feqn4}
g(s-t)=h(s)f(t),\ s,t\in I,
\end{equation}
for some function $h$. This is a Cauchy type functional equation (see
\cite{AD}). It follows that $g$ is exponential, thus so is $f$. If $I=\mathbb R$,
then $f$ is not a probability density function; thus, as we claimed before,
either $I=(-\infty,0)$ or $I=(0,+\infty)$.

\textbf{Case (ii)}: $|A|\ne|B|$. We first show that $I=\mathbb R$. Suppose indirectly that
$f(t_0)=0$ for some $t_0$. Then from \eqref{feqn2} we get
\[
f\!\big((A+B)s+(A-B)t_0\big)=0
\]
for all values of $s$, that is, $f\equiv 0$, which is impossible.

Now we put $f(x)=f(p(x))$, $q(x)=\log g(x)$, $C=\log 2A$ in \eqref{feqn1} to get
\begin{equation}\label{feqn5}
p\Big(\frac{u+v}{2}\Big)+p\Big(\frac{u-v}{2}\Big)=p(Au + Bv)+q(v)+C.
\end{equation}
Let $\Omega(u)$ a nonnegative, infinitely differentiable function such that
\[
\Omega(u)=0\mbox{ for }|u|\ge 1,\mbox{ and }\int_{-\infty}^{+\infty}\Omega(u)\,\mathrm du=1.
\]
Put $\Omega_{\varepsilon}(u)=\dfrac{1}{\varepsilon}\Omega\Big(\dfrac{u}{\varepsilon}\Big)$.
Multiplying both sides of \eqref{feqn5} by $\Omega_{\varepsilon}(x-u)$ and integrating with
respect to $u$ over the whole real line we get
\begin{equation}\label{feqn6}
p_{\varepsilon}\Big(\frac{x+v}{2}\Big)+p_{\varepsilon}\Big(\frac{x-v}{2}\Big)=
p_{\varepsilon}(Ax + Bv)+q(v)+C,
\end{equation}
where $p_{\varepsilon}$ is the convolution of $p$ and $\Omega_{\varepsilon}$. It is well known
that $p_{\varepsilon}$ is infinitely differentiable. Differentiate both sides of \eqref{feqn6}
with respect to $x$ to get
\[
p'_{\varepsilon}\Big(\frac{x+v}{2}\Big)+p'_{\varepsilon}\Big(\frac{x-v}{2}\Big)=
2A\,p'_{\varepsilon}(Ax + Bv).
\]
Set  $s=\dfrac{x+v}{2}$, $t=\dfrac{x-v}{2}$. Then
\[
p'_{\varepsilon}(s)+p'_{\varepsilon}(t)=2A\,p'_{\varepsilon}\big((A+B)s+(A-B)t\big).
\]
Differentiating both sides with respect to $s$ we get 
\[
p''_{\varepsilon}(s)=2A(A+B)\,p''_{\varepsilon}\big((A+B)s+(A-B)t\big).
\]
Here the left hand side does not depend on $t$. Therefore $p''_{\varepsilon}$ is constant and
$p_{\varepsilon}$ is a quadratic polynomial. Finally, $p_{\varepsilon}\to p$ a.e.\ as $\varepsilon\to 0$,
hence $p$ is also a quadratic polynomial a.e., thus $f=\exp(p)$ is normal.
\qed

\section{Independence of symmetrized random variables}

Let the random vectors $(X,Y)$ and $(X',Y')$ be independent and identically distributed. Suppose that 
the symmetrized $X-X'$ and $Y-Y'$ are independent. Does it follow that $X$ and $Y$ 
are also independent?

The answer is negative as it is shown by the following
 
\subsubsection*{Counterexample}

Let the random variables $X$ and $Y$ have joint characteristic function
\[
\hat f(t,s)= e^{-2|t|-t^2+it^2s}/(1+s^2).
\]
Later we will prove that this is in fact a characteristic function. It is easy to see that
\[
|\hat f(t,s)|=|\hat f(t,0)||\hat f(0,s)|
\]
for all $t,s$, thus if $(X',Y')$ is an independent copy of $(X,Y)$ then $X-X'$ and 
$Y-Y'$ are independent. On the other hand, $X$ and $Y$ are not independent because
\[
\hat f(t,s)\not\equiv\hat f(t,0)\hat f(0,s).
\]
What remains to be proved is that $\hat f(t,s)$ is in fact a characteristic function, or 
equivalently, the corresponding inverse Fourier transform
\[
f(x,y)=\frac{1}{(2\pi)^2}\int_{-\infty}^{\infty}\int_{-\infty}^{\infty}
\hat f(t,s)e^{-itx-isy}\,\mathrm dt\,\mathrm ds
\]
is nonnegative. Now
\[
f(x,y)=\frac{1}{4\pi}\int_{-\infty}^{\infty}e^{-2|t|-t^2-itx}
\int_{-\infty}^{\infty}\frac{e^{is(t^2-y)}}{\pi(1+s^2)}\,\mathrm ds\,\mathrm dt
\]
The second integral is the characteristic function of the standard Cauchy distribution
at $t^2-y$, thus
\begin{align*}
f(x,y)&=\frac{1}{4\pi}\int_{-\infty}^{\infty}e^{-2|t|-t^2-itx}\,
e^{-|t^2-y|}\,\mathrm dt\\[3pt]
&=\frac{1}{2}\,e^{-|y|}\,\frac{1}{2\pi}\int_{-\infty}^{\infty}
e^{-2|t|-t^2-|t^2-y|+|y|}\,e^{-itx}dt.
\end{align*}
It is enough to show that 
\[
\hat g(t)=e^{-2|t|-t^2-|t^2-y|+|y|}
\]
is a characteristic function for every real number $y$, because then the inverse Fourier 
transform
\[
g(x)=\frac{1}{2\pi}\int_{-\infty}^{\infty}\hat g(t)\,e^{-itx}\,\mathrm dt
\]
is nonnegative and so is $f(x,y)=\frac{1}{2}\,e^{-|y|}g(x)$. (Note that $\hat g(t)$
and hence $g(x)$ depend on $y$ as well.)

There are two cases. If $y\le 0$ then $\hat g(t)=e^{-2|t|-2t^2}$ is obviously a characteristic 
function, because it is the product of two characteristic functions.
The other case is where $y>0$. Then $\hat g(t)$ is a Polya type characteristic function 
meaning that $\hat g(0)=1$, $\hat g(-t)=\hat g(t)$, $g(t)\to 0$ non-increasingly
as $t\to\infty$ and $g(t)$ is a convex function on the positive real half-line. The proof
of all claims are easy. We just need to separate the cases where $0\le t\le\sqrt{y}$ 
and $t>\sqrt{y}$. In the former case $\hat g(t)=e^{-2t}$, while in 
the latter case $\hat g(t)=e^{-2t-2(t^2-y)}$. Clearly, $g'(t)\le 0$ 
and $g''(t)\ge 0$ for all $t>0$. 

\medskip

\begin{remark}
Let $X$ and $Y$ be real valued random variables with finite second moments.
The distance covariance of $X$ and $Y$ can be defined in the 
following form \cite{SzR,SzRB}.
 
Let $(X,Y)$, $(X',Y')$, $(X'',Y'')$ denote independent and identically distributed 
copies then the distance covariance is the square root of
\begin{align*}
{\mathop{\mathrm{dCov}}}^2(X,Y)&:= \E(|X -X'||Y -Y'|) + \E(|X-X'|
\E(|Y-Y'|)\\  
&\qquad -2 \E(|X-X'||Y-Y''|)\\
&=\cov(|X -X'|,|Y -Y'|)-2\cov(|X-X'|,|Y-Y''|).
\end{align*}
The definition of their \textit{distance correlation} is the following:
\[
\mathop{\mathrm{dCor}}(X,Y):=\frac{\mathop{\mathrm{dCov}}(X,Y)}
{\sqrt{\rule{0pt}{9pt}\mathop{\mathrm{dCov}}(X,X)\mathop{\mathrm{dCov}}(Y,Y)}}, 
\]
provided the denominator is positive.
An important property of distance correlation is that it characterizes independence, i.e. 
$\mathop{\mathrm{dCor}}(X,Y) = 0$ if and only if $X$ and $Y$ are independent.
The counterexample above makes the distance correlation even more interesting.
\end{remark}

\medskip

Next we show that the independence of $X$ and $Y$ does follow if in addition we
suppose that $X$ and $Y$ are bounded and their joint distribution is symmetric in the
sense that $(X, Y)$ and $(-X,-Y)$ are indentically distributed. What we really need
is less than boundedness; we only suppose that all (mixed) moments are finite and they
determine the joint distribution. 
\begin{theorem}
Let the random vectors $(X,Y)$ and $(X',Y')$ be identically distributed with bounded and 
symmetric joint distribution. Suppose that $X-X'$ and $Y-Y'$ are independent. 
Then so are $X$ and $Y$.
\end{theorem} 
\proof
Let us start from the Taylor expansion of the joint characteristic function,
\[
\hat f(t,s)=\E\left(e^{itX+isY}\right)=\sum_{k=0}^{\infty}\sum_{\ell=0}^{\infty}
\frac{(it)^k(is)^{\ell}}{k!\,\ell !}\E(X^k Y^{\ell}).
\]
It is enough to show that $\E(X^k Y^{\ell})=\E(X^k)\E(Y^{\ell})$ for all 
$k,\ell=1,2,\dots$, because then
\[
\hat f(t,s)=\sum_{k=0}^{\infty}\frac{(it)^k}{k!}\E(X^k)\cdot
\sum_{\ell=0}^{\infty}\frac{(is)^{\ell}}{\ell !}\E(Y^{\ell})=\hat f(t,0)\hat f(0,s),
\]
thus the independence of $X$ and $Y$ follows.

By symmetry, we have $\E(X^k)=0$ and $\E(Y^{\ell})=0$ if $k$ and $\ell$ are
odd. Similarly, $\E(X^k Y^{\ell})=0=\E(X^k)\E(Y^{\ell})$ if $k+\ell$ is odd. 
Hence all we have to show is that $\E(X^kY^{\ell})=\E(X^k)\E(Y^{\ell})$ 
whenever $k+\ell$ is even. This can be done by induction over $n:=k+\ell$.
Let $k=\ell=1$. Then
\[
\E\left((X-X')(Y-Y')\right)=\E(XY)-\E(X'Y)-\E(XY')+\E(YY')=2\E(XY),
\]
and by the independence of $X-X'$ and $Y-Y'$, 
\[
\E\left((X-X')(Y-Y')\right)=\E(X-X')\E(Y-Y')=0.
\]
Next, let $k+\ell=n$ even, and suppose that $\E(X^iY^j)=\E(X^i)\E(Y^j)$ if
$i+j<n$. Then
\begin{multline*}
\E\left((X-X')^k(Y-Y')^{\ell}\right)\\
\begin{aligned}
&=\sum_{i=0}^k\sum_{j=0}^{\ell}(-1)^{n-i-j}\binom{k}{i}\binom{\ell}{j}
\E\left(X^i(X')^{k-i}Y^j(Y')^{\ell-j}\right)\\
&=\sum_{i=0}^k\sum_{j=0}^{\ell}(-1)^{n-i-j}\binom{k}{i}\binom{\ell}{j}
\E\left(X^iY^j\right)\E\left((X')^{k-i}(Y')^{\ell-j}\right)
\end{aligned}
\end{multline*}
Now by the induction hypothesis we have
\begin{multline*}
\E\left((X-X')^k(Y-Y')^{\ell}\right)=
2\E(X^kY^{\ell})-2\E(X^k)\E(Y^{\ell})\\
+\sum_{i=0}^k\sum_{j=0}^{\ell}(-1)^{n-i-j}\binom{k}{i}\binom{\ell}{j}
\E(X^i)\E(X^{k-i})\E(Y^j)\E(Y^{\ell-j}).
\end{multline*}
Here the last line is equal to
\begin{gather*}
\sum_{i=0}^k(-1)^{k-i}\binom{k}{i}\E(X^i)\E(X^{k-i})\cdot
\sum_{j=0}^{\ell}(-1)^{\ell-j}\binom{\ell}{j}\E(Y^j)\E(Y^{\ell-j})\\
=\sum_{i=0}^k(-1)^{k-i}\binom{k}{i}\E\left(X^i(X')^{k-i}\right)\cdot
\sum_{j=0}^{\ell}(-1)^{\ell-j}\binom{\ell}{j}\E\left(Y^j(Y')^{\ell-j}\right)\\[3pt]
=\E\left((X-X')^k\right)\E\left((Y-Y')^{\ell}\right)\\[3pt]
=\E\left((X-X')^k(Y-Y')^{\ell}\right).
\end{gather*}
Therefore
\[
\E(X^kY^{\ell})=\E(X^k)\E(Y^{\ell})
\]
as needed.
\qed

\begin{remark}
A similar result is proved in a forthcoming paper by Jakob 
Raymaekers and Peter Rousseeuw \cite{RR}.
\end{remark}

\section{The Hungarian roulette}
\suppressfloats[t]

In this section we are going to pose an unsolved problem. Though the model to be 
introduced is quite simple, we couldn't find it in the literature.

$n$ people stand in a circle; everyone has a gun. At a given signal, everybody 
selects someone randomly and shoots at him. Each shot is fatal. The survivors 
continue this over and over again until there is one single survivor left or none. 
Let $p_n$ denote the probability that someone remains standing
in the end, and $\xi_n$ the number of people staying alive after the first round.
Then $p_0=0$, $p_1=1$, $p_2=0$, and 
\[
p_n=\sum_{k=0}^{n-2}p_k\,\P(\xi_n=k),\quad n\ge 3.
\]

\begin{figure}[h]
\centering
\includegraphics[scale=.6]{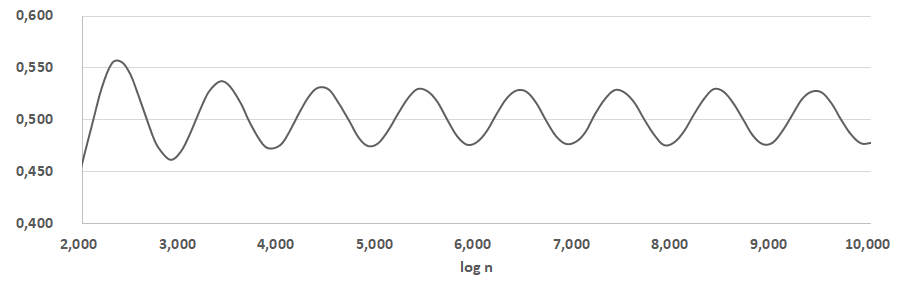}
\caption{Is this logarithmic periodicity or convergence?}
\label{Fig1}
\end{figure}

What is the asymptotic behavior of $p_n$ as $n\to\infty$? Approximately 
computing the values by Monte Carlo and plotting the probabilities against 
$\log n$ one can see sine-like periodicity, though with waves with slightly 
decreasing amplitude (see Figure \ref{Fig1}). Does it even subsist asymptotically 
or diminishes to convergence? We do not know the answer yet, but we will discuss 
the subject.
\smallskip

Fix $n$ and number the players from $1$ to $n$. Now, $\xi_n=\sum_{i=1}^n Y_i$,
where $Y_i=\ind\left(\mbox{player $i$ survives the first round}\right)$. 
These indicators are exchangeable and 
\[
\E\big(Y_{i_1}\cdots Y_{i_\ell}\big)=\left(\frac{n-\ell}{n-1}\right)^{\!\!\ell}
\left(\frac{n-\ell-1}{n-1}\right)^{\!\!n-\ell}
\]
for $1\le i_1<\dots<i_\ell\le n$. Therefore
\[
\P(\xi_n=k)=\sum_{\ell=k}^{n-2}(-1)^{\ell-k}\binom{\ell }{ k}
\binom{n }{ \ell}\frac{(n-\ell)^\ell(n-\ell-1)^{n-\ell}}{(n-1)^n}
\]
by Waring's formula (\cite{Feller}, Section IV.3). Moreover,
\[
\E(Y_i)=\left(\frac{n-2}{n-1}\right)^{\!\!n-1},\quad
\E(Y_i Y_j)=\left(\frac{n-2}{n-1}\right)^{\!\!2}\left(\frac{n-3}{n-1}\right)^{\!\!n-2},\ i\ne j,
\]

\begin{gather*}
\E(\xi_n/n)=\left(\frac{n-2}{n-1}\right)^{\!\!n-1},\\
\Var(\xi_n/n)=\frac{1}{n}\left(\frac{n-2}{n-1}\right)^{\!\!n-1}+
\frac{n-1}{n}\left(\frac{n-2}{n-1}\right)^{\!\!2}\left(\frac{n-3}{n-1}\right)^{\!\!n-2}-
\left(\frac{n-2}{n-1}\right)^{\!\!2n-2}.
\end{gather*}
By Taylor expansion it is easy to see that
\[
\left(\frac{n-a}{n-b}\right)^{\!\!n-c}=e^{-a+b}\left(1-\frac{(a-b)(a+b-2c)}{2n}+
O\big(n^{-2}\big)\right).
\]
Therefore we have
\[
\E(\xi_n/n)=\frac{1}{e}+O\big(n^{-1}\big),\quad
\Var(\xi_n/n)= \frac{1}{n}\left(\frac{1}{e}-\frac{2}{e^2}\right)+O\big(n^{-2}\big).
\]

Furhermore,
\begin{equation}\label{CLT}
\sqrt{n}\left(\frac{\xi_n}{n}-\frac{1}{e}\right)\to N\left(0,\,
\frac{1}{e}-\frac{2}{e^2}\right)
\end{equation}
in distribution. This can be proved as follows.

Let us slightly modify the game. Allow people to shoot at themselves, and let
$\xi'_n$ the number of players remainig alive. This is just the classical random 
allocation problem where $n$ balls are thrown into $n$ urns uniformly at random,
independently of each other. Here $\xi'_n$ is the number of empty urns. Then
by Theorem 1 of Ch.\ I, \S 3 in \cite{KSC} $\xi'_n$ satisfies \eqref{CLT}. Let 
us make a coupling between $\xi_n$ and $\xi'_n$ in the following way. Denote
by $\eta$ the number of people who turn the gun on themselves in the modified 
model. We ask them to shoot someone else instead (chosen from the others
with equal probability). In this way we just obtain the original model. Each
redirection changes the number of survivning players by at most one, thus
$|\xi_n-\xi'_n|\le\eta$. Since $\E(\eta)=1$, \eqref{CLT} immediately follows for 
$\xi_n$.

Denote how far counterclockwise the person player $i$ is targeting is by $X_i$ 
(thus $X_i\in\{1,2,\dots, n-1\}$). Then $X_1,\dots,X_n$ are iid, and
with a suitable function $\varphi{:}\{1,\dots,n-1\}^n\to\{0,1,\dots,n-2\}$ we
have $\xi_n=\varphi(X_1,\dots,X_n)$. Obviously, if one of the arguments
$x_1,\dots,x_n$ is changed (i.e., one of the players shoots at somebody else)
then $\varphi$ can change by at most one. Now the McDiarmid inequality
\cite{McDiarmid} gives
\[
\P\big(|\xi_n-\E(\xi_n)|\ge\varepsilon\big)\le 2\,\exp\big({-2}\varepsilon^2/n\big).
\]
This shows that $p_n$ is basically equal to the weighted mean of $o(n)$ 
probabilities $p_k$ around $p_{n/e}$, which makes it likely that $p_n$
shows some kind of logarithmic periodicity. Of course, this allows convergence, too.

\begin{conj}
The sequence $p_n$ does not converge as $n\to\infty$ over the positive integers, 
but it does along all subsequences of $(n)$ for which the fractional part of $\log n$ 
converges (asymptotic logarithmic periodicity).
\end{conj}
This phenomenon is called \emph{merging} \cite{DDF}.

In what follows we are going to present a heuristic reasoning to show why we refute 
convergence. We argue that the increasingly flattening waves do not flatten out completely.

Suppose $p_n\approx c+h(\log n)$, where $h$ is a smooth, sine-like damped periodic
 function with period $1$, and $c$ is a constant between the peaks and
the troughs, i.e., $h$ is not of constant sign. Consider a whole wave (period) over an 
interval of unit length. Let $x=M$ be the point where $f(x)$ is maximal in that interval 
and let $n=e^M$. Then
\[
h(M+1)\approx p_{en}-c=\E\left(p_{\xi_{en}}-c\right).
\]
By \eqref{CLT} we have
\[
\xi_{en}=n+\sigma\sqrt{n}\,\zeta_n=n\left(1+\frac{\sigma}{\sqrt{n}}\,\zeta_n\right),
\]
where $\sigma=\sqrt{1-\frac{2}{e}}$ and $\zeta_n\to N(0,1)$ in distribution.
Therefore
\[
\log\xi_{en}\approx \log n+\frac{\sigma}{\sqrt{n}}\,\zeta_n=
M+\sigma e^{-M/2}\zeta_n,
\]
thus
\[
h(M+1)\approx\E h\left(M+\sigma e^{-M/2}\,\zeta_n\right).
\]
Since $h'(M)=0$  we get
\[
h(M+1)\approx h(M)+\frac{\sigma^2}{2}e^{-M}f''(M)\E(\zeta^2_n)
\approx h(M)+\kappa f''(M)e^{-M},
\]
with a suitable constant $\kappa$. Note that $f''(M)\le 0$.

If the waves are similar in shape, i.e., each sufficiently distant wave is a
constant times of the previous one, then  $f''(x+k)/f(x+k)$, $k=0,1,\dots$ are 
approximately equal, thus
\[
h(M+1)\approx h(M)\left(1-\kappa'e^{-M}\right),\quad\kappa'\ge 0,
\]
which is a lower bound for the maximum of the next wave. Therefore, the
wave maximum $k$ waves away is no less than
\[
h(M+k)\approx h(M)\prod_{i=0}^{k-1}\left(1-\kappa'e^{-M-i}\right).
\]
Since $\sum_ie^{-M-i}<\infty$, these maxima do not drop down to $0$, not even 
gradually. It can be seen in exactly the same way that the minima do not increase 
to $0$ either.

\begin{remark}
Asymptotic logarithmic periodicity emerges in many different probabilistic models. 
An old example is the length $L_n$ of the longest head run in a sequence of $n$ 
tosses with a fair coin. It satisfies
\[
\P(L_n-\lfloor\log_2 n\rfloor<k)=\exp\left(-2^{-(k+1-\{\log_2 n\}}\right)+o(1),
\]
for every integer $k$, as $n\to \infty$, where $\{a\}=a-\lfloor a\rfloor$, the fractional 
part of $a$, see \cite{G}. A recent publication is \cite{HHKT}, where the discussed model 
also seems to exhibit log-periodicity, though the authors only give a heuristic proof. See 
also the references therein.
\end{remark}

\vspace{1cm}

\noindent\textbf{Tam\'as F.~M\'ori}\\
HUN-REN Alfr\'ed R\'enyi Institute of Mathematics\\
Budapest\\
Hungary\\
{\tt mori.tamas@renyi.hu}

\vspace{3ex}
\noindent\textbf{G\'abor J.~Sz\'ekely}\\
HUN-REN Alfr\'ed R\'enyi Institute of Mathematics\\
Budapest\\
Hungary\\
{\tt gabor.j.szekely@gmail.com}

\end{document}